\documentclass[12pt,a4paper]{amsart}

\def\frak{\mathfrak}

\newtheorem{theorem}{Theorem}[section]
\newtheorem{lemma}[theorem]{Lemma}
\newtheorem{corollary}[theorem]{Corollary}
\newtheorem{conjecture}[theorem]{Conjecture}

\theoremstyle{definition}

\theoremstyle{remark}
\newtheorem{remark}[theorem]{Remark}

\numberwithin{equation}{section}

\def\br{\mathbb R}
\def\bc{\mathbb C}

\def\mbx#1{\makebox[#1cm]{}}
\def\codim{\,\mbox{codim}\,}

\def\bp{\begin{pmatrix}}
\def\ep{\end{pmatrix}}
\def\op{{\rm op}}

\begin{document}

\title[Equivalence of domains II]{Equivalence of domains arising from duality of orbits on flag manifolds II}

\author{Toshihiko MATSUKI}
\address{Department of Mathematics\\
        Faculty of Science\\
        Kyoto University\\
        Kyoto 606-8502, Japan}
\email{matsuki@math.kyoto-u.ac.jp}
\date{}

\begin{abstract}
In \cite{GM1}, we defined a $G_\br$-$K_\bc$ invariant subset $C(S)$ of $G_\bc$ for each $K_\bc$-orbit $S$ on every flag manifold $G_\bc/P$ and conjectured that the connected component $C(S)_0$ of the identity will be equal to the Akhiezer-Gindikin domain $D$ if $S$ is of nonholomorphic type. This conjecture was proved for closed $S$ in \cite{WZ1,WZ2,FH,M8} and for open $S$ in \cite{M8}. In this paper, we prove the conjecture for all the other orbits when $G_\br$ is of non-Hermitian type.
\end{abstract}

\maketitle

\section{Introduction}

Let $G_\bc$ be a connected complex semisimple Lie group and $G_\br$ a connected real form of $G_\bc$. Let $K_\bc$ be the complexification in $G_\bc$ of a maximal compact subgroup $K$ of $G_\br$. Let $X=G_\bc/P$ be a flag manifold of $G_\bc$ where $P$ is an arbitrary parabolic subgroup of $G_\bc$. Then there exists a natural one-to-one correspondence between the set of $K_\bc$-orbits $S$ and the set of $G_\br$-orbits $S'$ on $X$ given by the condition:
\begin{equation}
S\leftrightarrow S'\Longleftrightarrow S\cap S'\mbox{ is non-empty and compact} \tag{1.1}
\end{equation}
(\cite{M4}). For each $K_\bc$-orbit $S$ we defined in \cite{GM1} a subset $C(S)$ of $G_\bc$ by
$$C(S)=\{x\in G_\bc\mid xS\cap S'\mbox{ is non-empty and compact}\}$$
where $S'$ is the $G_\br$-orbit on $X$ given by (1.1).

Akhiezer and Gindikin defined a domain $D/K_\bc$ in $G_\bc/K_\bc$ as follows (\cite{AG}). Let $\frak{g}_\br=\frak{k}\oplus\frak{m}$ denote the Cartan decomposition of $\frak{g}_\br={\rm Lie}(G_\br)$ with respect to $K$. Let $\frak{t}$ be a maximal abelian subspace in $i\frak{m}$. Put
$$\frak{t}^+=\{Y\in\frak{t}\mid |\alpha(Y)|<{\pi\over 2} \mbox{ for all }\alpha\in\Sigma\}$$
where $\Sigma$ is the restricted root system of ${\frak g}_\bc$
with respect to $\frak{t}$. Then $D$ is defined by
$$D=G_\br(\exp\frak{t}^+)K_\bc.$$

We conjectured the following in \cite{GM1}.

\begin{conjecture} {\rm (Conjecture 1.6 in \cite{GM1})} \ Suppose that $X=G_\bc/P$ is not $K_\bc$-homogeneous. Then we will have $C(S)_0=D$ for all $K_\bc$-orbits $S$ of nonholomorphic type on $X$. Here $C(S)_0$ is the connected component of $C(S)$ containing the identity. $($See \cite{GM1,M8} for the definition of the $K_\bc$-orbits of nonholomorphic type. When $G_\br$ is of non-Hermitian type, all the $K_\bc$-orbits are defined to be of nonholomorphic type.$)$
\end{conjecture}

Let $S_{\rm op}$ denote the unique open $K_\bc$-$B$ double coset in $G_\bc$ where $B$ is a Borel subgroup of $G_\bc$ contained in $P$. It is shown in \cite{H} and \cite{Ms} that $D\subset C(S_{\rm op})_0$. (The opposite inclusion $D\supset C(S_{\rm op})_0$ is proved in \cite{B}.) On the other hand the inclusion $C(S_{\rm op})_0\subset C(S)_0$ for every $K_\bc$-orbit $S$ on $X=G_\bc/P$ is shown in \cite{GM1} Proposition 8.1 and Proposition 8.3. So we have the inclusion
\begin{equation}
D\subset C(S)_0. \tag{1.2}
\end{equation}
We have only to prove the opposite inclusion. 

For a simple root $\alpha$ with respect to $B$ we can define a parabolic subgroup $P_\alpha$ by
$$P_\alpha=B\sqcup Bw_\alpha B$$
where $w_\alpha$ is the reflection for the root $\alpha$. Let $S_0$ be a closed $K_\bc$-$B$ double coset in $G_\bc$. Let $S_1,\ldots,S_\ell$ ($\ell=\codim_\bc S_0$) be a sequence of $K_\bc$-$B$ double cosets in $G_\bc$ and $\alpha_1,\ldots\alpha_\ell$ a sequence of simple roots such that
$$S_k^{cl}=S_0P_{\alpha_1}\cdots P_{\alpha_k}$$
and that
$$\dim_\bc S_k=\dim_\bc S_0+k$$
for $k=1,\ldots,\ell$ (c.f. \cite{GM2}, \cite{M3}, \cite{Sp}). Especially $S_\ell=S_\op$.

In this paper we first prove the following theorem.

\begin{theorem} \ Let $x$ be an element of $G_\bc$. If $I_0=xS_0\cap S'_\op P_{\alpha_\ell}\cdots P_{\alpha_1}$ is connected, then
$$I_k=xS_k^{cl}\cap S'_\op P_{\alpha_\ell}\cdots P_{\alpha_{k+1}}$$
is connected for $k=1,\ldots,\ell$. $($$S'_\op$ is the unique closed $G_\br$-$B$ double coset in $G_\bc$ which corresponds to $S_\op$ by $(1.1)$.$)$
\end{theorem}

\begin{remark} \ The sets $I_k\ (k=0,\ldots,\ell)$ are always nonempty because $xS_0P_{\alpha_1}\cdots P_{\alpha_\ell}=xS_\op^{cl}=G_\bc\supset S'_\op$.
\end{remark}

Let $S$ be a $K_\bc$-$P$ double coset in $G_\bc$. Then we can write
$$S^{cl}=S_k^{cl}=S_0P_{\alpha_1}\cdots P_{\alpha_k}$$
with some closed $K_\bc$-$B$ double coset $S_0$ and a sequence $\alpha_1,\ldots,\alpha_k$ of simple roots (\cite{M3}, \cite{Sp}). Secondly we prove the following.

\begin{theorem} \ {\rm (i)} \ If $x\in D^{cl}$, then $I_k$ is connected.

{\rm (ii)} \ If $x\in D^{cl}\cap C(S)$, then $I_k=xS\cap S'_k.$
\end{theorem}

As a corollary we solve Conjecture 1.1 for non-Hermitian cases:

\begin{corollary} \ Let $G_\br$ be simple and of non-Hermitian type. Then $C(S)_0=D$ for all the $K_\bc$-orbits $S\ne X$ on $X=G_\bc/P$.
\end{corollary}

\begin{proof} When $S$ is open in $G_\bc$, the equality is proved in \cite{M8}. So we may assume that $S$ is not open. Let $x$ be an element of $D^{cl}\cap C(S)$. Then we have only to show that $x\in D$.  Since $S_kP_{\alpha_{k+1}}\cdots P_{\alpha_{\ell-1}}\cap S_\op =\phi$, we have $S'_kP_{\alpha_{k+1}}\cdots P_{\alpha_{\ell-1}}\cap S'_\op =\phi$ by the duality (\cite{M2}) and therefore
$$S'_\op P_{\alpha_{\ell-1}}\cdots P_{\alpha_{k+1}}\cap S'_k=\phi.$$
By Theorem 1.4 (ii) we have
\begin{align*}
xS^{cl}\cap S'_\op P_{\alpha_{\ell-1}}\cdots P_{\alpha_{k+1}}
& = xS^{cl}\cap S'_\op P_{\alpha_\ell}\cdots P_{\alpha_{k+1}}\cap S'_\op P_{\alpha_{\ell-1}}\cdots P_{\alpha_{k+1}} \\
& = xS\cap S'_k\cap S'_\op P_{\alpha_{\ell-1}}\cdots P_{\alpha_{k+1}}= \phi.
\end{align*}
Hence
$$xS_{\ell-1}^{cl}\cap S'_\op =xS^{cl}P_{\alpha_{k+1}}\cdots P_{\alpha_{\ell-1}}\cap S'_\op =\phi.$$

For the orbit $S_{\ell-1}$ we defined the following domain $\Omega$ in \cite{GM2}.
$$\Omega=\{x\in G_\bc\mid xS_{\ell-1}^{cl}\cap S'_\op =\phi\}_0.$$
It is shown in \cite{FH} Theorem 5.2.6 and \cite{M8} Corollary 1.8 that
$$\Omega=D$$
when $G_\br$ is of non-Hermitian type. Hence $x\in D$.
\end{proof}

\begin{remark} \ Recently \cite{M10} proved Conjecture 1.1 for all non-closed $K_\bc$-orbits in Hermitian cases using Theorem 1.4. Thus the conjecture is now completely solved affirmatively.
\end{remark}

\section{$G_\br$-orbits on the full flag manifold}

The full flag manifold $\mathcal{F}$ of $G_\bc$ is the set of the Borel subgroups of $G_\bc$. If we take a Borel subgroup $B_0$ of $G_\bc$, then the factor space $G_\bc/B_0$ is identified with $\mathcal{F}$ by the map
$$G_\bc/B_0\ni gB_0\mapsto gB_0g^{-1}\in\mathcal{F}.$$
It is known that every $G_\br$-orbit ($G_\br$-conjugacy class) on $\mathcal{F}$ contains a Borel subgroup of the form
$$B=B(\frak{j},\Sigma^+)=\exp\left(\sum_{\alpha\in\Sigma^+\sqcup\{0\}} \frak{g}_\bc(\frak{j},\alpha)\right)$$
where $\frak{j}$ is a $\theta$-stable Cartan subalgebra of $\frak{g}_\br$, $\Sigma^+$ is a positive system of the root system $\Sigma$ of the pair $(\frak{g}_\bc, \frak{j}_\bc)$ and $\frak{g}_\bc(\frak{j},\alpha)=\{X\in\frak{g}_\bc\mid [Y,X]=\alpha(Y)X\mbox{ for all }Y\in\frak{j}\}$ (\cite{A}, \cite{M1}, \cite{R}).

Roots in $\Sigma$ are usually classified as follows.

(i) \ If $\theta(\alpha)=\alpha$ and $\frak{g}_\bc(\frak{j},\alpha)\subset \frak{k}_\bc$, then $\alpha$ is called a ``compact root''.

(ii) \ If $\theta(\alpha)=\alpha$ and $\frak{g}_\bc(\frak{j},\alpha)\subset \frak{m}_\bc$, then $\alpha$ is called a ``noncompact root''.

(iii) \ If $\theta(\alpha)=-\alpha$, then $\alpha$ is called a ``real root''.

(iv) \ If $\theta(\alpha)\ne \pm\alpha$, then $\alpha$ is called a ``complex root''.

For a simple root $\alpha$ of $\Sigma^+$ define the parabolic subgroup $P_\alpha$ as in Section 1. By the same arguments as in \cite{V} Lemma 5.1 and \cite{M3} Lemma 3, we can prove the following decomposition of $P_\alpha/B\cong P^1(\bc)$ into the $P_\alpha\cap G_\br$-orbits.

\begin{lemma} \ {\rm (i)} \ If $\alpha$ is compact, then $P_\alpha=(P_\alpha\cap G_\br)B$.

{\rm (ii)} \ If $\alpha$ is noncompact or real, then $P_\alpha/B\cong P^1(\bc)=\bc\sqcup\{\infty\}$ is decomposed into the three $(P_\alpha\cap G_\br)_0$-orbits $H_+,\ H_-$ and $H_0$ which are diffeomorphic to the upper half plane, the lower half plane and $P^1(\br)=\br\sqcup\{\infty\}$, respectively. $($Sometimes $H_+$ and $H_-$ are in the same $P_\alpha\cap G_\br$-orbit.$)$

{\rm (iii)} If $\alpha$ is complex, then $P_\alpha/B$ is decomposed into the two $P_\alpha\cap G_\br$-orbits consisting of a point $yB$ and the complement $(P_\alpha-yB)/B$.
\end{lemma}

\begin{remark} \ Concerning the $K_\bc$-action on $G_\bc/B$, it is shown in \cite{V} Lemma 5.1 (c.f. \cite{M3} Lemma 3, \cite{GM1} Lemma 9.1) that:

(i) \ If $\alpha$ is compact, then $P_\alpha=(P_\alpha\cap K_\bc)B$.

(ii) \ If $\alpha$ is noncompact or real, then $P_\alpha/B$ is decomposed into three $(P_\alpha\cap K_\bc)_0$-orbits consisting of two points and the complement.

(iii) If $\alpha$ is complex, then $P_\alpha/B$ is decomposed into two $P_\alpha\cap K_\bc$-orbits consisting of a point and the complement.
\end{remark}

As a corollary of Lemma 2.1 we have:

\begin{corollary} \ Let $g$ be an arbitrary element of $G_\bc$. Then every $(gP_\alpha g^{-1}\cap G_\br)_0$-invariant closed subset of $gP_\alpha/B$ is connected.
\end{corollary}

\begin{remark} \ On the contrary a $gP_\alpha g^{-1}\cap K_\bc$-invariant closed subset of $gP_\alpha/B$ may not be connected in view of Remark 2.2 (ii). 
\end{remark}

\section{Proof of the theorems}

\noindent {\it Proof of Theorem 1.2.} \ We will prove the theorem by induction on $k$. Suppose that $I_{k-1}$ is connected. Then
\begin{align*}
I_{k-1}P_{\alpha_k} & =(xS_{k-1}^{cl}\cap S'_\op P_{\alpha_\ell}\cdots P_{\alpha_k})P_{\alpha_k} \\
& =xS_k^{cl}\cap S'_\op P_{\alpha_\ell}\cdots P_{\alpha_k} \\
& =(xS_k^{cl}\cap S'_\op P_{\alpha_\ell}\cdots P_{\alpha_{k+1}}) P_{\alpha_k} \\
& =I_kP_{\alpha_k}
\end{align*}
is connected. Suppose that $I_k=A_1\sqcup A_2$ with some nonempty closed subsets $A_1$ and $A_2$ of $I_k$. Then we will get a contradiction. Since the Borel subgroup $B$ is connected, $A_1$ and $A_2$ are right $B$-invariant. Since $A_1P_{\alpha_k}$ and $A_2P_{\alpha_k}$ are closed and
$$A_1P_{\alpha_k}\cup A_2P_{\alpha_k}=I_kP_{\alpha_k}$$
is connected, we have $A_1P_{\alpha_k}\cap A_2P_{\alpha_k}\ne\phi$. Take an element $g$ of $A_1P_{\alpha_k}\cap A_2P_{\alpha_k}$. Then $gP_{\alpha_k}\cap I_k$ is decomposed as
$$gP_{\alpha_k}\cap I_k=(gP_{\alpha_k}\cap A_1)\sqcup (gP_{\alpha_k}\cap A_2)$$
with two nonempty closed subsets $gP_{\alpha_k}\cap A_1$ and $gP_{\alpha_k}\cap A_1$. But this contradicts Corollary 2.3 because $gP_{\alpha_k}\cap I_k=gP_{\alpha_k}\cap S'_\op P_{\alpha_\ell}\cdots P_{\alpha_{k+1}}$ is $gP_{\alpha_k}g^{-1}\cap G_\br$-invariant.

\hfill$\square$

\begin{lemma} \ {\rm (i)} \ $S_k$ is relatively closed in $S_\op P_{\alpha_\ell}\cdots P_{\alpha_{k+1}}$.

{\rm (ii)} \ $S'_k$ is relatively open in $S'_\op P_{\alpha_\ell}\cdots P_{\alpha_{k+1}}$.
\end{lemma}

\begin{proof} By the duality for the closure relation (\cite{M3}) we have only to show (i). Let $\widetilde{S}$ be a $K_\bc$-$B$ double coset contained in the boundary of $S_k$. Then
$$\codim_\bc \widetilde{S}>\codim_\bc S_k=\ell-k.$$
Hence $\widetilde{S}$ cannot be contained in $S_\op P_{\alpha_\ell}\cdots P_{\alpha_{k+1}}$ by \cite{V} Lemma 5.1 (c.f. \cite{GM1} Lemma 9.1).
\end{proof}

\noindent {\it Proof of Theorem 1.4.} \ (i) \ Since $S_0/B$ is compact and $S'_0/B$ is open, we see that
\begin{align*}
C(S_0) &=\{x\in G_\bc\mid xS_0\cap S'_0\mbox{ is nonempty and closed in }G_\bc\} \\
&=\{x\in G_\bc\mid xS_0\subset S'_0\}.
\end{align*}
Hence $C(S_0)_0$ is the cycle space for $S'_0$ defined in \cite{WW}. Since $D\subset C(S_0)_0$ by (1.2), it follows that
$$x\in D\Longrightarrow xS_0\subset S'_0.$$

Suppose that $x\in D^{cl}$. Then we have
$$xS_0\subset {S'_0}^{cl}\subset S'_\op P_{\alpha_\ell}\cdots P_{\alpha_1}$$
and hence $I_0=xS_0$ is connected. By Theorem 1.2 the intersection
$$I_k=xS^{cl}\cap S'_\op P_{\alpha_\ell}\cdots P_{\alpha_{k+1}}$$
is connected.

(ii) \ By Lemma 3.1 $S'_k$ is relatively open in $S'_\op P_{\alpha_\ell}\cdots P_{\alpha_{k+1}}$. On the other hand $S$ is also relatively open in $S^{cl}$. Hence
\begin{equation}
\mbox{$xS\cap S'_k$ is relatively open in $I_k=xS^{cl}\cap S'_\op P_{\alpha_\ell}\cdots P_{\alpha_{k+1}}$.} \tag{3.1}
\end{equation}

Suppose that $x\in C(S)$. Then $xS\cap S'$ is nonempty and closed in $G_\bc$ by definition. Since $S'_k$ is relatively closed in $S'$, it follows that
\begin{equation}
xS\cap S'_k\mbox{ is closed in }G_\bc. \tag{3.2}
\end{equation}
Since $xS\cap S'=(xS\cap S'_k)P$, it also follows that
\begin{equation}
xS\cap S'_k\mbox{ is nonempty}. \tag{3.3}
\end{equation}

Suppose moreover that $x\in D^{cl}$. Then $I_k$ is connected by (i). Hence it follows from (3.1), (3.2) and (3.3) that
$$\mbx6 I_k=xS\cap S'_k.\mbx6\square$$


\begin{thebibliography}{BGW}

\bibitem[A]{A}
K. Aomoto,
\textit{On some double coset decompositions of complex semi-simple Lie groups},
J. Math. Soc. Japan \textbf{18} (1966), 1--44.

\bibitem[AG]{AG}
D. N. Akhiezer and S. G. Gindikin,
\textit{On Stein extensions of real symmetric spaces},
Math. Ann. \textbf{286} (1990), 1--12.

\bibitem[B]{B}
L. Barchini,
\textit{Stein extensions of real symmetric spaces and the geometry of the flag manifold},
Math. Ann. \textbf{326} (2003), 331--346.

\bibitem[FH]{FH}
G. Fels and A. Huckleberry,
\textit{Characterization of cycle domains via Kobayashi hyperbolicity},
preprint (AG/0204341).

\bibitem[GM1]{GM1}
S. Gindikin and T. Matsuki,
\textit{Stein extensions of Riemannian symmetric spaces and dualities of orbits on flag manifolds},
Transform. Groups \textbf{8} (2003), 333--376.

\bibitem[GM2]{GM2}
S. Gindikin and T. Matsuki,
\textit{A remark on Schubert cells and duality of orbits on flag manifolds},
preprint (RT/0208071).

\bibitem[H]{H}
A. Huckleberry,
\textit{On certain domains in cycle spaces of flag manifolds},
Math. Ann. \textbf{323} (2002), 797--810.

\bibitem[M1]{M1}
T. Matsuki,
\textit{The orbits of affine symmetric spaces under the action of minimal parabolic subgroups},
J. Math. Soc. Japan \textbf{31} (1979), 331--357.

\bibitem[M2]{M2}
T. Matsuki,
\textit{Orbits on affine symmetric spaces under the action of parabolic subgroups},
Hiroshima Math. J. \textbf{12} (1982), 307--320.

\bibitem[M3]{M3}
T. Matsuki,
\textit{Closure relations for orbits on affine symmetric spaces under the action of minimal parabolic subgroups},
Adv. Stud. Pure Math. \textbf{14} (1988), 541--559.

\bibitem[M4]{M4}
T. Matsuki,
\textit{Closure relations for orbits on affine symmetric spaces under the action of parabolic subgroups. {I}ntersections of associated orbits},
Hiroshima Math. J. \textbf{18} (1988), 59--67.

\bibitem[M5]{Ms}
T. Matsuki,
\textit{Stein extensions of Riemann symmetric spaces and some generalization},
J. Lie Theory \textbf{13} (2003), 563--570.

\bibitem[M6]{M8}
T. Matsuki,
\textit{Equivalence of domains arising from duality of orbits on flag manifolds},
preprint (RT/0309314).

\bibitem[M7]{M10}
T. Matsuki,
\textit{Equivalence of domains arising from duality of orbits on flag manifolds III},
preprint (RT/0410302).

\bibitem[R]{R}
W. Rossmann,
\textit{The structure of semisimple symmetric spaces},
Canad. J. Math. \textbf{31} (1979), 157--180.

\bibitem[Sp]{Sp}
T. A. Springer,
\textit{Some results on algebraic groups with involutions},
Adv. Stud. Pure Math. \textbf{6} (1984), 525--534.

\bibitem[V]{V}
D. A. Vogan,
\textit{Irreducible characters of semisimple Lie groups} III,
Invent. Math. \textbf{71} (1983), 381--417.

\bibitem[WW]{WW}
R. O. Wells and J. A. Wolf,
\textit{Poincar{\'e} series and automorphic cohomology on flag domains},
Ann. of Math. \textbf{105} (1977), 397--448.

\bibitem[WZ1]{WZ1}
J. A. Wolf and R. Zierau,
\textit{Linear cycle spaces in flag domains},
Math. Ann. \textbf{316} (2000), 529--545.

\bibitem[WZ2]{WZ2}
J. A. Wolf and R. Zierau,
\textit{A note on the linear cycle spaces for groups of Hermitian type},
J. Lie Theory \textbf{13} (2003), 189--191.

\end{thebibliography}
\end{document}